\documentclass[11pt]{article}

\usepackage{amsfonts}
\usepackage{amsmath}
\usepackage{enumerate}
\usepackage[parfill]{parskip}
\usepackage{color}

\newtheorem{theorem}{Theorem}

\newtheorem{example}{Example}

\newtheorem{lemma}{Lemma}

\newtheorem{proposition}{Proposition}

\newenvironment{proof}[1][Proof]{\textbf{#1.} }{\ \rule{0.5em}{0.5em}}
\input{tcilatex}
\begin{document}

\title{Discretely Sampled Variance and Volatility Swaps versus their
Continuous Approximations}
\author{Robert Jarrow\thanks{%
Johnson Graduate School of Management, Cornell University, Ithaca, NY, 14853
and Kamakura Corporation} \and Younes Kchia\thanks{%
Centre de Math\'ematiques Appliqu\'ees, Ecole Polytechnique, Paris} \and %
Martin Larsson\thanks{%
School of Operations Research, Cornell University, Ithaca, NY, 14853} \and %
Philip Protter\thanks{%
Statistics Department, Columbia University, New York, NY, 10027} \thanks{%
Supported in part by NSF grant DMS-0906995}}
\date{\today}
\maketitle

\begin{abstract}
Discretely sampled variance and volatility swaps trade actively in OTC
markets. To price these swaps, the continuously sampled approximation is
often used to simplify the computations. The purpose of this paper is to
study the conditions under which this approximation is valid. Our first set
of theorems characterize the conditions under which the discretely sampled
swap values are finite, given the values of the continuous approximations
exist. Surprisingly, for some otherwise reasonable price processes, the
discretely sampled swap prices do not exist, thereby invalidating the
approximation. Examples are provided. Assuming further that both swap values
exist, we study sufficient conditions under which the discretely sampled
values converge to their continuous counterparts. Because of its popularity
in the literature, we apply our theorems to the 3/2 stochastic volatility
model. Although we can show finiteness of all swap values, we can prove
convergence of the approximation only for some parameter values.

KEY\ WORDS: variance swaps, volatility swaps, NFLVR, semimartingales
\end{abstract}

\section{Introduction}

Although variance and volatility swaps only started trading in the
mid-1900s, they have since become a standard financial instrument useful for
managing volatility risk (see Carr and Lee \cite{Carr/Lee:2009} for the
history of volatility derivative markets). In the pricing and hedging of
variance and volatility swaps, a distinction is made between payoffs that
are discretely or continuously sampled. Discretely sampled variance and
volatility swaps trade in the over-the counter (OTC) markets. In contrast,
continuously sampled variance and volatility swaps are only an abstract
construct, often used to approximate the values of their discretely sampled
counterparts (see Broadie and Jain \cite{Broadie/Jain:2008},\cite%
{Broadie/Jain 2008a}, Chan and Platen \cite{ChanPlaten 2009}, and Carr and
Lee \cite{Carr/Lee: 2009a}). These approximations depend on an exchange of a
limit operator (as the discrete sampling period goes to zero) with an
expectation operator. This operator exchange is often invoked without
adequate justification. The purpose of this paper is to characterize the
conditions under which this operator exchange is valid. 

For our investigation we utilize the martingale pricing methodology where we
take as given the asset's price process assuming markets are arbitrage free
(in the sense of No Free Lunch with Vanishing Risk (NFLVR)). This evolution
is taken to be very general. We only assume that the price process is a
strictly positive semimartingale with possibly discontinuous sample paths.
We also assume, as a standing hypothesis, that the continuously sampled
variance and volatility swaps have finite values. Otherwise, before the
analysis begins, the approximation would not make sense.

Our first two Theorems (1 and 2) characterize the \textit{additional conditions} needed on the price process such that the discretely sampled variance and volatility swaps have finite values. Of course, when the discretely sampled variance and volatility swap values do not exist, the approximation is again nonsensical. Surprisingly, we provide examples of otherwise reasonable price processes, with stochastic volatility of the volatility, where these discretely sample variance swap values do not exist.

Next, assuming both the continuous and discretely sampled variance swap values are finite, we study conditions justifying an exchange of the limit and expectation operators. In this regard, under no additional hypotheses, Theorem 3 provides an upper bound for the maximum difference between these two values. Theorem 4 proves, under some additional moment conditions, that the exchange of the two operators is valid. Furthermore, this theorem also provides infomation on the rate of convergence ($1/n$ where $n$ is the number of discretely sampled prices).

Lastly, given the recent interest in the 3/2 stochastic volatility model for pricing volatility derivatives (see Carr and Sun \cite{Carr/Sun:2007}, Chan and Platen \cite{ChanPlaten 2009}), we explore its consistency with valuing discretely sampled variance swaps with their continuously
sampled counterparts. Here we show that both the discrete and continuously sampled variance swaps have finite values. Unfortunately, we can only prove convergence of the discrete to the continuously sampled variance swap values for some parameter ranges, but not all. A complete characterization of this convergence for the 3/2 stochastic volatility model remains an open question.

An outline for this paper is as follows. Section 2 gives the framework underlying the model. Section 3 studies the finiteness and convergence of
the discretely sampled variance and volatility swap values. Finally, Section 4 provides examples to illustrate the theorems proved.

\section{Framework}

Let a filtered probability space $(\Omega ,\mathcal{F},\mathbb{F},P)$ be
given, where the filtration $\mathbb{F}=(\mathcal{F}_{t})_{t\in \lbrack 0,T]}
$ satisfies the usual conditions and $T$ is a fixed time horizon. We suppose
that there is a liquidly traded asset paying no dividends, whose market
price process is modeled by a semimartingale $S=(S_{t})_{t\in \lbrack 0,T]}$
such that $S>0$ and $S_->0$. The value of the money market account is chosen
as numeraire, or, viewed differently, the interest rate is zero. The price
process is assumed to be arbitrage free in the sense of No Free Lunch with
Vanishing Risk (NFLVR), see Delbaen and Schachermayer~\cite%
{Delbaen/Schachermayer:1994} and~\cite{Delbaen/Schachermayer:1998}, which
guarantees the existence of at least one equivalent probability measure
under which $S$ is a local martingale. We assume $P$ is such a measure, and
that the modeler prices future payoffs by taking expectations with respect
to~$P$.

\subsection{Variance and volatility swaps}

A \textit{variance swap} (with strike zero) and maturity $T$ is a contract
which pays the "realized variance," i.e. the square of the logarithm returns
up to time $T$, namely 
\begin{equation*}
\frac{252}{n}\sum_{i=1}^{n}\left( \ln \frac{S_{t_{i}}}{S_{t_{i-1}}}\right)
^{2},
\end{equation*}%
where $0=t_{0}<...<t_{n}=T$ is a regular sampling of the time interval $%
\left[ 0,T\right] $, i.e.~$t_i-t_{i-1}=\frac{1}{n}$ for $i=1,\ldots,n$.
Finally, $252$ is the number of trading days per year. The maturity $T$ is
approximately $\frac{n}{252}$.

Rather than considering the quantity 
\begin{equation*}
P^{n}(T)=\sum_{i=1}^{n}\left( \ln \frac{S_{t_{i}}}{S_{t_{i-1}}}\right) ^{2},
\end{equation*}%
practitioners\footnote{%
Confirmed in a private discussion with Peter Carr.} often use its limit 
\begin{equation*}
P(T)=[\ln S,\ln S]_{T}
\end{equation*}%
in the pricing of variance swaps (see Carr and Lee~\cite{Carr/Lee:2009}).
That is, one approximates the quantity $E(P^n(T))$ by $E(P(T))$. The
continuous approximation $P(T)$ to the variance swap's true payoff $P^n(T)$
is justified by the fact that 
\begin{equation*}
[\ln S,\ln S]_{T}=\lim_{n\rightarrow \infty }\sum_{i=1}^{n}\left( \ln \frac{%
S_{t_{i}}}{S_{t_{i-1}}}\right) ^{2},
\end{equation*}%
where the limit is in probability and taken over a sequence of subdivisions
whose mesh size tends to zero.

Analogous to a variance swap, a \textit{volatility swap} is a security
written on the square root of the variance swap's payoff. We will use the
notation 
\begin{equation*}
V^{n}(T)=\sqrt{\sum_{i=0}^{n-1}\left( \ln \frac{S_{t_{i}}}{S_{t_{i-1}}}%
\right) ^{2}}\qquad \text{and}\qquad V(T)=\sqrt{\langle \ln S\rangle _{T}}.
\end{equation*}%
Again, $V^{n}(T)$ is the volatility swap payoff up to multiplicative and
additive constants.

Our aim is to investigate the validity of these approximations, motivated by
the fact that the convergence of $P^{n}(T)$ to $P(T)$, respectively $V^n(T)$
to $V(T)$, is only in probability, which \emph{a priori} does not guarantee
convergence of their expectations. It is the latter convergence one needs in
order to justify the use of this approximation in the context of pricing.

To simplify the notation later on, we introduce the following convention.

\textbf{Notation.} For a process $X=(X_t)_{t\in[0,T]}$ we define 
\begin{equation*}
\delta_i X = X_{t_i}-X_{t_{i-1}} \qquad (i=1,\ldots,n). 
\end{equation*}
In particular, $P^n(T)=\sum_{i=1}^{n} (\delta_i \ln S)^2$.

\subsection{Mathematical preliminaries}

\label{S:prel}

Since $S>0$ and $S_->0$ there is a semimartingale $M$ such that 
\begin{equation*}
S = S_0 \mathcal{E}(M), 
\end{equation*}
where $\mathcal{E}(\cdot)$ denotes stochastic exponential; see~\cite%
{Jacod/Shiryaev:2003}, Theorem~II.8.3. In fact, since $S$ is a local
martingale and hence a special semimartingale, it follows that from
Theorem~II.8.21 in~\cite{Jacod/Shiryaev:2003} that $M$ is a stochastic
integral with respect to the local martingale part of $S$. Furthermore,
since $S>0$, the jumps of $M$ satisfy $\Delta M > -1$, so by the
Ansel-Stricker Theorem (see \cite{Ansel/Stricker:1994}), $M$ is a local
martingale.

The following notation is standard. We direct the reader to~\cite%
{Jacod/Shiryaev:2003} for details. The random measure $\mu^M$ associated
with the jumps of $M$ is given by 
\begin{equation*}
\mu^M(dt,dx) = \sum_s {\boldsymbol{1}_{\{\Delta M_s \neq 0\}}}%
\varepsilon_{(s,\Delta M_s)}(dt,dx), 
\end{equation*}
and its predictable compensator is $\nu(dt,dx)$. We may then write 
\begin{equation}  \label{eq:S}
S = S_0 \exp\left\{ M - \frac{1}{2}\langle M^c,M^c\rangle - (x-\ln(1+x)) *
\mu^M \right\},
\end{equation}
and furthermore decompose $M$ as 
\begin{equation}  \label{eq:M}
M = M^c + M^d = M^c + x*(\mu^M-\nu),
\end{equation}
where $M^c$ is the continuous local martingale part of $M$ and $%
M^d=x*(\mu^M-\nu)$ is the jump part. Both these processes are local
martingales. The following elementary result will be useful later, so we
state it as a lemma.

\begin{lemma}
\label{L:QV} The quadratic variation of $M$ and $\ln S$ are given by

\begin{itemize}
\item[(i)] $[M,M] = \langle M^c,M^c\rangle + x^2 * \mu^M$

\item[(ii)] $[\ln S,\ln S] = \langle M^c, M^c \rangle + (\ln(1+x))^2 * \mu^M$
\end{itemize}
\end{lemma}

\begin{proof}
Part $(i)$ is a standard fact. For part $(ii)$, write $X=(x-\ln(1+x)) * \mu^M
$ and notice that $\ln S - \ln S_0 = M-X- \frac{1}{2}\langle M^c,M^c\rangle$%
. Hence 
\begin{align*}
[\ln S,\ln S] &= [M-X,M-X] \\
&= [M,M] + [X,X] - 2[M,X] \\
&= \langle M^c,M^c\rangle + x^2 * \mu^M + (x-\ln(1+x))^2 * \mu^M - 2
x(x-\ln(1+x)) * \mu^M.
\end{align*}
By the Cauchy-Schwartz inequality, 
\begin{align*}
|x(x-\ln(1+x))| * \mu^M &= \sum_{s\leq \cdot} \left|\Delta M_s ( \Delta M_s
- \ln(1+\Delta M_s)) \right| \\
& \leq \sqrt{ \sum_{s\leq \cdot} (\Delta M_s)^2 } \sqrt{ \sum_{s\leq \cdot}
( \Delta M_s - \ln(1+\Delta M_s))^2} <\infty,
\end{align*}
so $x(x-\ln(1+x)) * \mu^M$ also converges absolutely, a.s. Termwise
manipulation is therefore allowed, so since $x^2 + (x-\ln(1+x))^2 -
2x(x-\ln(1+x)) = (\ln(1+x))^2$, we get 
\begin{equation*}
x^2*\mu^M + (x-\ln(1+x))^2 * \mu^M - 2 x(x-\ln(1+x)) * \mu^M = (\ln(1+x))^2
* \mu^M. 
\end{equation*}
The result follows.
\end{proof}

We also have the following lemma, which gives the semimartingale
decomposition of~$\ln S$.

\begin{lemma}
\label{L:vf} Assume that $\ln S$ is locally integrable. Then 
\begin{equation*}
\ln S -\ln S_0 = M^c + \ln(1+x)*(\mu^M-\nu) - \frac{1}{2}\langle
M^c,M^c\rangle - (x-\ln(1+x)) * \nu. 
\end{equation*}
\end{lemma}

\begin{proof}
From (\ref{eq:S}) and (\ref{eq:M}) we have 
\begin{equation*}
\ln S -\ln S_0 = M^c - \frac{1}{2}\langle M^c,M^c\rangle + x*(\mu^M-\nu) -
(x-\ln(1+x)) * \mu^M. 
\end{equation*}
Our assumption together with the fact that $M^c$ and $x*(\mu^M-\nu)$ are
local martingales, hence locally integrable, implies that $\langle
M^c,M^c\rangle$ and $(x-\ln(1+x)) * \mu^M$ are locally integrable (notice
that both are nonnegative). Hence $(x - \ln(1+x)) * \nu$ is locally
integrable (see~\cite{Jacod/Shiryaev:2003}, Proposition~II.1.8), so we may
add and subtract this quantity to the right side of the previous display to
obtain 
\begin{equation*}
\ln S -\ln S_0 = M^c - \frac{1}{2}\langle M^c,M^c\rangle +
\ln(1+x)*(\mu^M-\nu) - (x-\ln(1+x)) * \nu. 
\end{equation*}
This is the desired expression.
\end{proof}

\section{Approximation using the quadratic variation}

\subsection{Finiteness of expectations}

\label{S:fe}

In order for there to be any hope that $E(P(T))$ accurately approximates $%
E(P^n(T))$, a minimal requirement is that both these quantities be finite.
Perhaps somewhat surprisingly, they need not be finite or infinite
simultaneously. This is of course a potentially serious issue, since the
value of $E(P(T))$, when finite, is nonsensical as an approximation of $%
E(P^n(T))$ if this is infinite.

The following result gives necessary and sufficient conditions for $E(P^n(T))
$ to be finite, given that the approximation $P(T)$ is known to have finite
expectation.

\begin{theorem}
\label{T:int1} Assume that $P(T)\in L^1$. The following statements are
equivalent.

\begin{itemize}
\item[(i)] $P^n(T) \in L^1$ for at least one $n\geq 1$

\item[(ii)] $\left\{ 
\begin{array}{l}
\langle M^c,M^c\rangle_T \in L^2 \\ 
\! \\ 
(x - \ln(1+x)) * \nu_T \in L^2%
\end{array}%
\right.$
\end{itemize}
\end{theorem}

\begin{proof}
First of all, note that since $P^n(T)=\sum_{i=1}^n (\delta_i \ln S)^2$, we
have $P^n(T)\in L^1$ if and only if $\delta_i \ln S\in L^2$ for each $i$,
which is equivalent to having $\ln S_{t_i} - \ln S_0\in L^2$ for each~$i$.
We thus need to show that this is equivalent to condition~$(ii)$ in the
statement of the theorem.

By Lemma~\ref{L:QV}, our basic assumption $P(T)=[\ln S,\ln S]_T\in L^1$
implies that $\langle M^c,M^c\rangle_T\in L^1$ and $(\ln(1+x))^2*\mu^M_T\in
L^1$. Hence both $M^c_t$ and $\ln(1+x)*(\mu^M-\nu)_t$ are in $L^2$ for every 
$t\leq T$. Therefore, using the representation from Lemma~\ref{L:vf}, namely 
\begin{equation*}
\ln S -\ln S_0 = M^c - \frac{1}{2}\langle M^c,M^c\rangle +
\ln(1+x)*(\mu^M-\nu) - (x-\ln(1+x)) * \nu, 
\end{equation*}
we deduce that $\ln S_{t_i}-\ln S_0\in L^2$ for each $i$ holds if and only
if 
\begin{equation*}
\frac{1}{2}\langle M^c,M^c\rangle_{t_i} + (x-\ln(1+x)) * \nu_{t_i}\in L^2 
\end{equation*}
for each $i$. Since both $\langle M^c,M^c\rangle$ and $(x-\ln(1+x)) * \nu$
are nonnegative and nondecreasing, this is equivalent to condition~$(ii)$.
The proof is finished.
\end{proof}

We note that the above theorem implies that if $P^n(T) \in L^1$
for some $n$ then $P^n(T) \in L^1$ for all $n\geq 1$. An analogous result
holds for volatility swaps.

\begin{theorem}
\label{T:int1vol} Assume that $V(T)\in L^1$. The following statements are
equivalent.

\begin{itemize}
\item[(i)] $V^n(T) \in L^1$ for some $n\geq 1$

\item[(ii)] $\left\{ 
\begin{array}{l}
\langle M^c,M^c\rangle_T \in L^1 \\ 
\! \\ 
(x - \ln(1+x)) * \nu_T \in L^1%
\end{array}%
\right.$
\end{itemize}
\end{theorem}

\begin{proof}
Using for instance the fact that all norms on $\mathbb{R}^n$ are equivalent,
there are constants $0<c\leq C<\infty$ such that 
\begin{equation*}
cV^n(T)\leq\sum_{i=1}^{n}\big|\delta_i \ln S\big|\leq C V^n(T) \quad \text{%
a.s.} 
\end{equation*}
Thus $V^n(T)\in L^1$ is equivalent to $\delta_i \ln S \in L^1$ for each $i$,
which is equivalent to $\ln S_{t_i} - \ln S_0\in L^1$ for each~$i$.

Now the proof is similar to that of Theorem~\ref{T:int1}. By Lemma~\ref{L:QV}%
, $V(T)=\sqrt{[\ln S,\ln S]_T} \in L^1$ implies that $\sqrt{\langle
M^c,M^c\rangle_T}$ and $\sqrt{(\ln(1+x))^2*\mu^M_T}$ are in $L^1$, which via
the Burkholder-Davis-Gundy inequalities implies that $M^c_t$ and $\ln(1+x) *
(\mu^M-\nu)_t$ are in $L^1$ for each $t\leq T$. By Lemma~\ref{L:vf},
therefore, $\ln S_{t_i} - \ln S_0\in L^1$ for each~$i$ if and only if 
\begin{equation*}
\frac{1}{2}\langle M^c,M^c\rangle_{t_i} + (x-\ln(1+x)) * \nu_{t_i}\in L^1 
\end{equation*}
for each $i$. This is equivalent to condition~$(ii)$ of the theorem.
\end{proof}

Again, we note that Theorem~\ref{T:int1vol} implies that if $%
V^n(T) \in L^1$ for some $n$ then $V^n(T) \in L^1$ for all $n\geq 1$.


When $S$ is continuous we have the following variation of Theorem~\ref{T:int1}, which indicates that in
many cases, $E(P^n(T))<\infty$ is a stronger requirement than $E(P(T))<\infty
$. This will be used in Section~\ref{S:ex}, where we discuss specific examples. Notice that when $S$ is continuous, $P(T)=[\ln S,\ln S]=\langle
M^c,M^c\rangle = \langle M,M\rangle$.

\begin{proposition}
\label{P:fp} Assume that $S$ is continuous. Then the following are
equivalent.

\begin{itemize}
\item[(i)] $M_{t}\in L^{1}$ and $P^{n}(T)\in L^{1}$

\item[(ii)] $P(T)\in L^{2}$
\end{itemize}
\end{proposition}

\begin{proof}
In the continuous case, $\ln S - \ln S_0 = M - \frac{1}{2}\langle M,M\rangle$%
. First assume~$(ii)$, i.e.~that $P(T)=\langle M,M\rangle_T\in L^2$. Then
certainly $M_t\in L^2$ for every $t\leq T$, hence $\ln S_t\in L^2$ for every 
$t\leq T$. Therefore $\delta_i \ln S\in L^2$ for every $i$, and as in the
beginning of the proof of Theorem~\ref{T:int1}, this implies that $P^n(T)\in
L^1$. Hence $(i)$ holds. Conversely, if~$(i)$ is satisfied, then $%
\delta_i\ln S\in L^2$ for every $i$, so $\ln S_T \in L^2$. We also have $%
M_T\in L^1$ by assumption, so $\langle M,M\rangle_T=2(M_T - \ln S_T + \ln
S_0)$ is in $L^1$, implying that $M_T$ is in fact in $L^2$. This lets us
strengthen the previous conclusion to $\langle M,M\rangle_T\in L^2$, which
is~$(ii)$.
\end{proof}

Notice that the implication $(ii)\Rightarrow(i)$ also follows from Theorem~%
\ref{T:int1}. We can prove an analogous result for volatility swaps.

\begin{proposition}
Assume that $S$ is continuous. Then the following are equivalent.

\begin{itemize}
\item[(i)] $M_{t}\in L^{1}$ and $V^{n}(T)\in L^{1}$

\item[(ii)] $V(T)\in L^{2}$, i.e. $\langle M,M\rangle_T \in L^1$.
\end{itemize}
\end{proposition}

\begin{proof}
As in the proof of Theorem~\ref{T:int1vol}, $V^n(T)\in L^1$ is equivalent to 
$\delta_i \ln S \in L^1$ for each $i$. The rest of the proof is similar to
that of Proposition~\ref{P:fp}.
\end{proof}

In the case of continuous price processes, Proposition~\ref{P:fp} makes it
clear how one can proceed to construct examples where the approximation $P(T)
$ has finite expectation, but the true payoff $P^n(T)$ does not. Indeed, any
process $S$ of the form $S=S_0\mathcal{E}(M)$ will do, where $M$ is a
continuous local martingale that satisfies 
\begin{equation}  \label{eq:f}
\left\{ 
\begin{array}{l}
\langle M,M\rangle_T \in L^1 \\ 
\! \\ 
\langle M,M\rangle_T \notin L^2.%
\end{array}%
\right.
\end{equation}

It is clear that such processes exist; what is less clear is to what extent
examples can be found among models that appear in applications. In Section~%
\ref{S:ex} we provide examples demonstrating that the condition~(\ref{eq:f})
can appear in models that may appear innocuous at first sight.

\subsection{Bounds on the approximation error}

\label{S:bds}

In this section we assume that both $E(P(T))$ and $E(P^n(T))$ are finite,
and study conditions under which they are close for large $n$; as already
mentioned, although $P^n(T) \to P(T)$ in probability, the expectations need
not converge. We start by showing that under general conditions, the two
expectations at least cannot be too far apart. We then impose additional
structure on the model and give conditions that guarantee convergence. The
focus of this section is on variance swaps. The analysis of volatility swaps is more complicated and is omitted.


\begin{theorem}
\label{T:bd1} Assume that $P^n(T)$ and $P(T)$ both are in $L^1$. Then there
is a constant $C>0$, independent of~$n$, such that 
\begin{equation*}
|E(P(T))-E(P^{n}(T))|\leq C\qquad \text{for\ all\ }n.
\end{equation*}
\end{theorem}

The proof of Theorem~\ref{T:bd1} is straightforward once the following lemma
has been established.

\begin{lemma}
\label{L:ci} Assume that $P^n(T)$ and $P(T)$ both are in $L^1$ and define 
\begin{align*}
N &= \ln(1+x) * (\mu^M-\nu) \\
A &= \frac{1}{2}\langle M^c,M^c\rangle + (x-\ln(1+x)) * \nu.
\end{align*}
Then $\langle M^c,M^c\rangle_T$ and $[N,N]_T$ are in $L^1$. Moreover, 
\begin{align*}
\big| E\big( P(T) \big) - E\big( P^n(T) \big) \big| &\leq E\Big( %
\sum_{i=1}^n (\delta_i A)^2 \Big) + 2 E\Big( \sum_{i=1}^n |\delta_iM^c|
\delta_i A\Big) + 2 E\Big(\sum_{i=1}^n |\delta_i N| \delta_i A \Big),
\end{align*}
and we have 
\begin{equation*}
E\Big( \sum_{i=1}^n |\delta_iM^c| \delta_i A\Big) \leq \sqrt{ E \big(\langle
M^c,M^c\rangle_T\big)}\sqrt{ E \big(\sum_{i=1}^n (\delta_i A)^2\big) } 
\end{equation*}
and 
\begin{equation*}
E\Big(\sum_{i=1}^n |\delta_i N| \delta_i A \Big) \leq \sqrt{ E \big(\lbrack
N,N]_T\big)}\sqrt{ E \big(\sum_{i=1}^n (\delta_i A)^2\big) }. 
\end{equation*}
\end{lemma}

\begin{proof}
Observe that $N$ is a purely discontinuous local martingale
with $[N,N]=(\ln(1+x))^2*\mu^M$, and $A$ is a nondecreasing process. Since
by assumption $P(T)=[\ln S,\ln S]_T \in L^1$, Lemma~\ref{L:QV} implies that
both $M^c$ and $N$ are $L^2$~martingales, which is the first assertion.
Since also $P^n(T)\in L^1$, Theorem~\ref{T:int1} shows that $A_T\in L^2$.
Again by Lemma~\ref{L:QV}, 
\begin{align*}
E\Big( \lbrack \ln S,\ln S]_T \Big) &= \sum_{i=1}^n E\Big( \delta_i \langle
M^c,M^c\rangle \Big) + \sum_{i=1}^n E\Big( \delta_i [N,N] \Big) \\
&= \sum_{i=1}^n E\Big( (\delta_i M^c)^2+ (\delta_i N)^2 \Big).
\end{align*}
Lemma~\ref{L:vf} shows that with $N$ and $A$ as above, we have $\ln S -\ln
S_0 = M^c + N - A$. Hence $\delta_i \ln S = \delta_i M^c + \delta_i N -
\delta_i A$, and therefore 
\begin{align*}
E\big( P^n(T) \big) - E\big( P(T) \big) &= \sum_{i=1}^n E\Big( (\delta_i
A)^2+ 2(\delta_i M^c) (\delta_i N) - 2 (\delta_iM^c) (\delta_i A) - 2
(\delta_i N)(\delta_i A) \Big) \\
&= \sum_{i=1}^n E\Big( (\delta_i A)^2 - 2 (\delta_iM^c) (\delta_i A) - 2
(\delta_i N)(\delta_i A) \Big),
\end{align*}
where the second equality holds because $M^c$ and $N$ are orthogonal $L^2$%
~martingales. The triangle inequality and Jensen's inequality yield 
\begin{align*}
\big| E\big( P(T) \big) - E\big( P^n(T) \big) \big| &\leq E\Big( %
\sum_{i=1}^n (\delta_i A)^2 \Big) + 2 E\Big( \sum_{i=1}^n |\delta_iM^c|
\delta_i A\Big) + 2 E\Big(\sum_{i=1}^n |\delta_i N| \delta_i A \Big) \\
&= (\text{I}) + (\text{II}) + (\text{III}).
\end{align*}
This settles the first inequality in the statement of the lemma. For (II),
applying the Cauchy-Schwartz inequality twice, first for the sum and then
the expectation, yields 
\begin{equation*}
(\text{II}) \leq 2 E\Bigg( \sqrt{ \sum_{i=1}^n (\delta_iM^c)^2 } \sqrt{
\sum_{i=1}^n (\delta_i A)^2 }\Bigg)
\leq 2 \sqrt{ E\Big( \sum_{i=1}^n (\delta_iM^c)^2\Big) } \sqrt{ E\Big( %
\sum_{i=1}^n (\delta_i A)^2 \Big) }. 
\end{equation*}
It only remains to notice that $E\Big( \sum_{i=1}^n (\delta_iM^c)^2\Big) =
\sum_{i=1}^n E\Big(\delta_i\langle M^c, M^c\rangle\Big) = E \big(\langle
M^c,M^c\rangle_T\big)$. An analogous calculation yields 
\begin{equation*}
(\text{III}) \leq 2 \sqrt{ E \big(\lbrack N,N]_T\big)}\sqrt{ E\Big( %
\sum_{i=1}^n (\delta_i A)^2 \Big) }, 
\end{equation*}
thus proving the lemma.
\end{proof}

\begin{proof}[Proof of Theorem~\protect\ref{T:bd1}]
By Lemma~\ref{L:ci} it suffices to bound $E\big( \sum_{i=1}^n (\delta_i A)^2 %
\big)$. Since $\delta_i A\geq 0$ for each $i$, we get 
\begin{equation*}
E\Big( \sum_{i=1}^n (\delta_i A)^2 \Big) \leq E\Big( \big\{ \sum_{i=1}^n
\delta_i A \big\}^2 \Big) = E(A_T^2) < \infty. 
\end{equation*}
The proof is complete.
\end{proof}

Notice that constant $C$ can be taken to be 
\begin{equation*}
C = E(A_T^2) + 2\sqrt{ E \big(\langle M^c,M^c\rangle_T\big)}\sqrt{ E(A_T^2)}
+ 2 \sqrt{ E \big(\lbrack N,N]_T\big)}\sqrt{ E(A_T^2)}, 
\end{equation*}
where $N$ and $A$ are as in Lemma~\ref{L:ci}. It is worth pointing out that
since $C$ is independent of $n$, it is in particular valid for $n=1$. It can
therefore not be expected to provide a tight bound for large~$n$. Such
results can be obtained under additional structure on the model, which we
now introduce.

\emph{For the remainder of this section we assume that our probability space supports an $m$-dimensional Brownian motion $W=(W^1,\ldots,W^n)$ and a
Poisson random measure $\mu = \mu(dt,dz)$ on $\mathbb{R}_+ \times \mathbb{R}$ with intensity measure $dt\otimes F(dz)$, where $\int (z^2\wedge 1)
F(dz)<\infty$. Moreover, we assume that $M^c$ is a stochastic integral with respect to $W$, and that $M^d$ is a stochastic integral with respect to $\mu
- dt\otimes F(dz)$. That is, we assume that there are predictable processes $a^1,\ldots,a^m$ and a predictable function $\psi>-1$ such that 
\begin{equation*}
M^c_t = \sum_{k=1}^m \int_0^t a^k_s dW^k_s \qquad \text{and} \qquad M^d_t =
\psi * (\mu - dt \otimes F(dz))_t. 
\end{equation*}
In this case the compensator $\nu$ of $\mu^M$ satisfies 
\begin{equation*}
\int_0^T \int_{(-1,\infty)} G(s,x) \nu(ds,dx) = \int_0^T \int_{\mathbb{R}}
G(s,\psi(s,z)) F(dz)ds 
\end{equation*}
for every nonnegative predictable function $G$. Moreover, it is a classical
result that there exists a one-dimensional Brownian motion $B$ and a
nonnegative predictable process $\sigma$ such that 
\begin{equation*}
M^c_t = \int_0^t \sigma_s dB_s. 
\end{equation*}%
}

Under this structure we can formulate conditions on $\sigma$ and $\psi$
under which the expectation of $P^n(T)$ converges to the expectation of $P(T)
$ as $n\to\infty$.

\begin{theorem}\label{T:c}
With the notation and assumptions just described
above, assume that the following conditions hold: 
\begin{equation*}
\left\{ 
\begin{array}{l}
E\Big\{ \int_0^T \sigma_s^4 ds \Big\} < \infty \\ 
\! \\ 
E\Big\{ \int_0^T \int_{\mathbb{R}} (1\vee|z|^{-p})\big\{ \psi(s,z)^2 +
(\ln(1+\psi(s,z)))^2 \big\} F(dz)ds \Big\} < \infty,%
\end{array}%
\right. 
\end{equation*}
for some $p\geq 0$ such that $\int_{\mathbb{R}} (1\wedge |z|^p) F(dz) <
\infty$. Then 
\begin{equation*}
\lim_{n\to\infty} | E(P(T)) - E(P^n(T)) | = 0. 
\end{equation*}
More specifically, there are constants $C$ and $D$ such that 
\begin{equation*}
| E(P(T)) - E(P^n(T)) | \leq Cn^{-1} + Dn^{-1/2}. 
\end{equation*}
\end{theorem}

Notice that $p\leq 2$ always works, in the sense that $\int_{\mathbb{R}}
(1\wedge z^2) F(dz) < \infty$. However, depending on $F(dx)$, smaller values
of $p$ may also work, imposing less stringent restrictions on~$\psi$. In
particular, if the Poisson random measure only has finitely many jumps, so
that $\int_{\mathbb{R}}F(dx)<\infty$, we may take $p=0$, and the condition
on $\psi$ reduces to 
\begin{equation*}
E\Big\{ \int_0^T \int_{\mathbb{R}} \big\{ \psi(s,z)^2 + (\ln(1+\psi(s,z)))^2 %
\big\} F(dz)ds \Big\} < \infty. 
\end{equation*}

We remark that Theorem~\ref{T:c} generalizes results in \cite%
{Broadie/Jain:2008}, where the authors study the Black-Scholes model, the
Heston stochastic volatility model, the Merton jump-diffusion model, and the
stochastic volatility with jumps model by Bates and Scott.

The proof of Theorem~\ref{T:c} requires two lemmas.

\begin{lemma}
\label{L:b1} For $0\leq s < t \leq T$ we have, a.s., 
\begin{align*}
\Big\{ \int_s^t \int_{(-1,\infty)} & (x - \ln(1+x)) \nu(du,dx) \Big\}^2 \\
&\leq (t-s) C \int_0^T \int_{\mathbb{R}} (1\vee|z|^{-p})\big\{ \psi(s,z)^2 +
(\ln(1+\psi(s,z)))^2 \big\} F(dz)ds
\end{align*}
for any $p\geq 0$ such that $\int_{\mathbb{R}} (1\wedge |z|^p) F(dz) < \infty
$, and a finite constant $C$ that does not depend on $s$, $t$ or $p$.
\end{lemma}

\begin{proof}
Jensen's inequality yields 
\begin{align*}
\Big\{ \int_s^t \int_{(-1,\infty)} (x - \ln & (1+x)) \nu(du,dx) \Big\}^2 \\
&= \Big\{ \int_s^t \int_{\mathbb{R}} (\psi(s,z) - \ln(1+\psi(s,z))) F(dz)ds %
\Big\}^2 \\
&\leq (t-s) \int_s^t \Big\{ \int_{\mathbb{R}} G(s,z) F(dz) \Big\}^2 ds,
\end{align*}
where we have defined $G(s,z)=\psi(s,z) - \ln(1+\psi(s,z))$. Splitting up
the integral over~$\mathbb{R}$ as the sum of the integrals over $\{|z|\leq
1\}$ and $\{|z|>1\}$, and applying the inequality $(a+b)^2\leq 2a^2 + 2b^2$,
we obtain 
\begin{equation*}
\Big\{ \int_{\mathbb{R}} G(s,z) F(dz) \Big\}^2 \leq 2\Big\{ \int_{\{|z|\leq
1\}} G(s,z) F(dz) \Big\}^2 + 2\Big\{ \int_{\{|z| > 1\}} G(s,z) F(dz) \Big\}%
^2. 
\end{equation*}
Since $C_1 = F(\{|z| > 1\})<\infty$, Jensen's inequality applied to the
second term yields 
\begin{equation*}
\Big\{ \int_{\{|z| > 1\}} G(s,z) F(dz) \Big\}^2 \leq C_1 \int_{\{|z| > 1\}}
G(s,z)^2 F(dz). 
\end{equation*}
For the first term, let $\widetilde F(dz) = |z|^p F(dz)$. Then $C_2 =
\widetilde F(\{|z| \leq 1\})<\infty$ by assumption, and using once again
Jensen's inequality we get 
\begin{align*}
\Big\{ \int_{\{|z|\leq 1\}} G(s,z) F(dz) \Big\}^2 &= \Big\{ \int_{\{|z|\leq
1\}} |z|^{-p} G(s,z) \widetilde F(dz) \Big\}^2 \\
&\leq C_2 \int_{\{|z| \leq 1\}} |z|^{-2p} G(s,z)^2 \widetilde F(dz) \\
&= C_2 \int_{\{|z| \leq 1\}} |z|^{-p} G(s,z)^2 F(dz).
\end{align*}
Assembling the different terms and noting that 
\begin{equation*}
G(s,z)^2 \leq 2\Big\{ \psi(s,z)^2 + (\ln(1+\psi(s,z)))^2 \Big\}
\end{equation*}
yields the result with $C=4(C_1\vee C_2)$.
\end{proof}

\begin{lemma}
\label{L:b2} Under the assumptions of Theorem~\ref{T:c}, we have $P(T)\in L^1
$, $\langle M^c,M^c\rangle_T\in L^2$, and $(x-\ln(1+x)) * \nu_T \in L^2$.
\end{lemma}

\begin{proof}
By Jensen's inequality, $E\{ (\int_0^T \sigma_s^2ds )^2 \} \leq
TE\{\int_0^T\sigma^4_sds\}$, which is finite by hypothesis. So $\langle
M^c,M^c\rangle_T\in L^2$. To prove that $P(T)\in L^1$ it therefore suffices
to note that 
\begin{equation*}
E\Big\{ (\ln(1+x))^2*\mu^M_T \Big\} = E\Big\{ (\ln(1+\psi))^2*\mu_T \Big\} =
E\Big\{ (\ln(1+\psi))^2*(ds\otimes F(dz)) \Big\}, 
\end{equation*}
which is finite by assumption, since $1\vee|z|^{-p} \geq 1$. Finally, an
application of Lemma~\ref{L:b1} with $s=0$ and $t=T$ gives a bound on $E\{ (
(x-\ln(1+x)) * \nu_T)^2\}$ that is finite by assumption. The lemma is proved.
\end{proof}

\begin{proof}[Proof of Theorem~\protect\ref{T:c}]
By Lemma~\ref{L:b2}, the hypotheses of Lemma~\ref{L:ci} are satisfied, so we
have the bound 
\begin{align*}
\big| E\big( P(T) \big) - E\big( P^n(T) \big) \big| &\leq E\Big( %
\sum_{i=1}^n (\delta_i A)^2 \Big) + 2 E\Big( \sum_{i=1}^n |\delta_iM^c|
\delta_i A\Big) + 2 E\Big(\sum_{i=1}^n |\delta_i N| \delta_i A \Big) \\
&= (\text{I}) + (\text{II}) + (\text{III}),
\end{align*}
where $N = \ln(1+x) * (\mu^M-\nu)$ and $A = \frac{1}{2}\langle
M^c,M^c\rangle + (x-\ln(1+x)) * \nu$. We first deal with~(I). Using the
inequality $(x+y)^2\leq 2x^2 + 2y^2$ we obtain 
\begin{equation*}
\sum_{i=1}^nE\Big( (\delta_i A)^2 \Big) \leq \sum_{i=1}^n \left\{ \frac{1}{2}
E\Big( (\delta_i \langle M^c,M^c\rangle)^2 \Big) + 2 E\Big( (\delta_i
(x-\ln(1+x))*\nu)^2 \Big) \right\}. 
\end{equation*}
Now, Jensen's inequality yields 
\begin{align*}
E\Big\{ (\delta_i \langle M^c,M^c\rangle)^2 \Big\} &= E\Big\{ \Big( %
\int_{t_{i-1}}^{t_i} \sigma_s^2ds \Big)^2 \Big\} \leq \frac{T}{n}E\Big\{ %
\int_{t_{i-1}}^{t_i} \sigma_s^4 ds \Big\}.
\end{align*}
Furthermore, Lemma~\ref{L:b1} with $s=t_{i-1}$ and $t=t_i$ yields 
\begin{align*}
E\Big\{ (\delta_i (x- & \ln(1+x))*\nu)^2 \Big\} \\
&\leq \frac{C_1}{n} \int_{t_{i-1}}^{t_i} \int_{\mathbb{R}} (1\vee|z|^{-p})%
\big\{ \psi(s,z)^2 + (\ln(1+\psi(s,z)))^2 \big\} F(dz)ds
\end{align*}
for some constant $C_1$ independent of $i$. Summing over $i$ shows that 
\begin{align*}
(\text{I}) \leq \frac{T}{2n} E\Big\{ & \int_0^T \sigma_s^4 ds \Big\} \\
&+ \frac{2C_1}{n}E\Big\{ \int_0^T \int_{\mathbb{R}} (1\vee|z|^{-p})\big\{ %
\psi(s,z)^2 + (\ln(1+\psi(s,z)))^2 \big\} F(dz)ds \Big\},
\end{align*}
which is equal to $n^{-1}$ times a constant $C$ that is finite by assumption.

Concerning the two remaining terms $(\text{II})$ and $(\text{III})$, Lemma~%
\ref{L:ci} gives 
\begin{equation*}
(\text{II}) \leq 2\sqrt{ E \big(\langle M^c,M^c\rangle_T\big)}\sqrt{ E \big(%
\sum_{i=1}^n (\delta_i A)^2\big) } 
\end{equation*}
and 
\begin{equation*}
(\text{III}) \leq 2 \sqrt{ E \big(\lbrack N,N]_T\big)}\sqrt{ E \big(%
\sum_{i=1}^n (\delta_i A)^2\big) }, 
\end{equation*}
so that $(\text{II})+(\text{III}) \leq C_2 \sqrt{ (\text{I}) } \leq C_2\sqrt{%
C} n^{-1/2}$ for some constant $C_2$. The claim now follows with $C$ as
above and $D=C_2\sqrt{C}$.
\end{proof}

\section{Examples}

\label{S:ex}

\subsection{Strict local martingales}

The literature on asset price bubbles centers around the phenomenon that $S$
can be a strict local martingale under the risk neutral measure~$P$, see~%
\cite{Cox/Hobson:2005}, \cite{Jarrow/Protter/Shimbo:2006},~\cite%
{Jarrow/Protter/Shimbo:2010}. Moreover, in~\cite{Ekstrom/Tysk:2009} this
issue has been noted to cause complications for pricing using PDE
techniques. On the other hand, alternative criteria of no arbitrage type
have been proposed by various authors to guarantee the existence of a \emph{%
true} martingale measure, for instance~\cite{Sin:1998}, \cite%
{Carr/Cherny/Urusov:2007},~\cite{Hobson:2010}. It is therefore natural to
ask about the relationship between our previous results and the true
martingale (or strict local martingale) property of $S$. We give two
examples, in the continuous case, showing that the two are not connected in
general. More specifically, the examples show that the martingale property
of $S$ has little to do with the integrability of $\langle M,M\rangle _{T}$.

Our first example uses the following criterion, which is well-known~\cite%
{Carr/Cherny/Urusov:2007}. As in Section~\ref{S:bds}, $B$ is standard
Brownian motion.

\begin{lemma}
\label{L:elv} Assume $S>0$ satisfies the stochastic differential equation $%
dS_t=\sigma(S_t)dB_t$. It is then a true martingale if and only if for some $%
a>0$, 
\begin{equation}
\int_{a}^{\infty }\frac{x}{\sigma ^{2}(x)}dx=\infty .  \label{eq:tm}
\end{equation}
\end{lemma}

\begin{example}[$S$ a strict local martingale and $\langle M,M\rangle_{T}$ $%
\in $ $L^{2}$]
Consider the Constant Elasticity of Variance (CEV) models $dS_{t}=S_{t}^{\alpha }dB_{t}$. By Lemma~\ref{L:elv}%
, $S$ is a strict local martingale if and only if $\alpha >1$. We would like
to choose $\alpha >1$ such that $M_{t}=\int_{0}^{t}S_{s}^{\alpha -1}dB_{s}$
is a martingale with an integrable quadratic variation, i.e.~$%
E\{(\int_{0}^{t}S_{s}^{2(\alpha -1)})^{2}\}<\infty $. This can be achieved
with $\varepsilon \in (0,1)$ and $\alpha =1+\frac{\varepsilon }{4}>1$.
Indeed, 
\begin{equation*}
E\left\{ \Big(\int_{0}^{T}S_{s}^{2(\alpha -1)}ds\Big)^{2}\right\} \leq
TE\left( \int_{0}^{T}S_{s}^{4(\alpha -1)}ds\right) =TE\left(
\int_{0}^{T}S_{s}^{\varepsilon }ds\right) =T\int_{0}^{T}E(S_{s}^{\varepsilon
})ds.
\end{equation*}%
Since $\varepsilon \in (0,1)$, $x\mapsto x^{\varepsilon }$ is concave.
Jensen's inequality thus implies that the right side above is dominated by $%
T\int_{0}^{T}E(S_{s})^{\varepsilon }ds\leq T^{2}S_{0}^{\varepsilon }<\infty $%
, where $E(S_s)\leq S_0$ because $S$ is a positive local martingale, hence a
supermartingale. This shows that $S$ can be a strict local martingale, even
if the quadratic variation $\langle M,M\rangle _{T}=\langle \ln S, \ln
S\rangle _{T}$ is in $L^{2}$.
\end{example}

For completeness, we also give a simple example showing that the reverse
situation is also possible: that $S$ can be well-behaved (a bounded
martingale), while $M$ is not.

\begin{example}[$S$ a bounded martingale and $M$ a strict local martingale]
Let $X$ be the reciprocal of a Bessel(3) process. It is well-known that $X$
is a strict local martingale, see e.g.~\cite[p. 20-21]{CW}. Set $M_{t}=-X_{t}
$, so that $M$ is a strict local martingale with $M_t\in L^1$ and $M_{t}<0$
a.s.~for all $t\geq 0$. Now, $P(T)=\langle M,M\rangle_T$ is not in $L^1$
(otherwise $M$ would be a true martingale), and also $P^n(T)$ fails to be in 
$L^1$ by Proposition~\ref{P:fp}.

However, since $S_{t}=\mathcal{E}(M)_t=\exp \{-X_{t}-\frac{1}{2}\langle
X,X\rangle _{t}\}\leq 1$, it is a bounded local martingale, hence a true
martingale. In this example, the \textquotedblleft bad\textquotedblright\
behavior of $M$ is caused by its ability to take on very large negative
values. This does not carry over to $S$, since it is obtained through
exponentiation.
\end{example}

\subsection{Stochastic volatility of volatility}

\label{S:svv}

We now proceed to give an example of a class of continuous stock price
models that look innocuous, but where the conditions~(\ref{eq:f}) at the end
of Section~\ref{S:fe} are satisfied for certain parameter values. In those
cases, Proposition~\ref{P:fp} implies that $E(P(T))<\infty$ but $%
E(P^n(T))=\infty$.

We use stochastic volatility models with stochastic volatility of
volatility. Let $B$, $W$ and $Z$ be three Brownian motions, and let $\rho$
denote the correlation between $W$ and~$Z$, i.e.~$d\langle W,Z\rangle_t=
\rho dt$. No restrictions will be imposed on the correlation structure of $%
(B,W,Z)$, other than through the parameter~$\rho$. We consider the following
model for the stock price $S$, its volatility $v$ and the volatility of
volatility~$w$: 
\begin{align}
dS_t &= S_t\sqrt{v_t}dB_t  \label{eph1} \\
dv_t &= v_t\sqrt{w_t}dW_t  \label{eph2} \\
dw_t &= \kappa (\theta - w_t)dt+\eta\sqrt{w_t}dZ_t ,  \label{eph3}
\end{align}
where $\kappa$, $\theta$, $\eta$ are positive constants. Maintaining our
previous notation, where $M$ is the stochastic logarithm of $S$, we have
that $M_t=\int_0^t\sqrt{v_s}dB_s$, so $\langle M, M\rangle_t=\int_0^t v_sds$%
. Recall condition~(\ref{eq:f}) from Section~\ref{S:fe}: 
\begin{equation*}
\left\{ 
\begin{array}{l}
\langle M,M\rangle_T \in L^1 \\ 
\! \\ 
\langle M,M\rangle_T \notin L^2.%
\end{array}%
\right. 
\end{equation*}

Note that $v$ is a nonnegative local martingale and hence a supermartingale.
Together with Fubini's theorem, this yields 
\begin{equation*}
E(\langle M, M \rangle_T) = E\Big\{\int_0^Tv_sds\Big\}=\int_0^TE(v_s)ds \leq
T v_0, 
\end{equation*}
so $\langle M, M\rangle_T \in L^1$. Now we wish to find conditions such that 
$\langle M, M\rangle_T= \int_0^T v_sds \notin L^2$. To this end, define 
\begin{equation*}
T^*=\sup\{t : E(v^2_t) < \infty\}, 
\end{equation*}
and let $\chi = 2\rho \eta-\kappa$ and $\Delta=\chi^2-2\eta^2$. It is proven
in \cite{Andersen/Piterbarg:2007} that 
\begin{equation*}
T^* = \left\{ 
\begin{array}{ll}
\frac{1}{\sqrt{\Delta}}\ln(\frac{\chi+\sqrt{\Delta}}{\chi-\sqrt{\Delta}}) & 
\text{if\ } \Delta\geq 0 \text{\ and\ } \chi>0, \\ 
\! &  \\ 
\frac{2}{\sqrt{-\Delta}}\Big(\arctan(\frac{\sqrt{-\Delta}}{\chi})+\pi 
\boldsymbol{1}_{\{\chi<0\}}\Big) & \text{if\ } \Delta < 0, \\ 
\! &  \\ 
+\infty & \text{otherwise}.%
\end{array}%
\right. 
\end{equation*}

The next step is to establish that $v$ is in fact a martingale.

\begin{lemma}
\label{L:vmg} The process $v$ defined above is a true martingale.
\end{lemma}

\begin{proof}
Using Feller's test of explosion (see e.g.~\cite{Rogers/Williams:1994}), a
straightforward calculation shows that $w$ does not explode under $P$.
Therefore, using the same techniques as in~\cite{Andersen/Piterbarg:2007},
it suffices to establish that the auxiliary process $\hat w$, defined as the
solution to 
\begin{equation*}
d\hat w_t=\big(\kappa(\theta-\hat w_t)+\rho\eta\hat w_t\big)dt+\eta\sqrt{%
\hat w_t}dZ_t, \qquad \hat w_0=w_0, 
\end{equation*}
is non-explosive. This can again be verified using Feller's criterion.
\end{proof}

We may now conclude our construction by choosing the parameters $\rho$, $\eta
$ and $\kappa$ such that $T^*<\infty$, and then choose $T>T^*$. In this
case, Fubini's theorem implies that 
\begin{equation*}
E\Big\{ \Big(\int_0^Tv_s ds\Big)^2\Big\}=E\Big\{\int_0^T\int_0^Tv_s v_t dsdt%
\Big\}=\int_0^T\int_0^TE(v_sv_t)dsdt. 
\end{equation*}
Moreover, $E(v_s v_t)=E(v_s E(v_t \mid \mathcal{F}_s))=E(v_s^2)$, so we get 
\begin{equation*}
E\Big\{ \Big(\int_0^Tv_s ds\Big)^2\Big\}=\int_0^T\int_0^T E(v_s^2) dsdt \geq
\int_{T^*}^T\int_{T^*}^T E(v_s^2)dsdt = \infty. 
\end{equation*}
We conclude that $\langle M, M\rangle_T \notin L^2$, and can summarize our
findings as follows.

\begin{example}
\label{ex:11} Suppose in the stock price model with stochastic volatility of
volatility described above, the parameters are such that $T>T^*$. Then the
preceding discussion shows that 
\begin{equation*}
E(P(T)) < \infty \qquad \text{but} \qquad E(P^n(T)) = \infty. 
\end{equation*}
That is, the approximation to the variance swap payoff has finite
expectation, whereas the true payoff does not.
\end{example}

It is interesting to note that it is sometimes possible to change to an
equivalent measure $Q\sim P$, under which the price process is still a local
martingale, and such that both $P(T)$ and $P^n(T)$ become integrable. The
authors would like to thank Kerry Back~\cite{KB} for posing the question of
whether or not this can happen. To carry out the construction, let us
continue to consider the stochastic volatility of volatility model described
above.

\begin{proposition}
Assume that we are in the framework of a doubly stochastic
volatility model as described in~\eqref{eph1},~\eqref{eph2}, and~\eqref{eph3}. Suppose that $\Delta\geq 0$ and $\chi>0$, so that $T^*<\infty$, and
assume also that $B$ is independent of $(W,Z)$. Then there is an equivalent
measure $Q$ such that $S$ is a local martingale under $Q$, and $\tilde
T^*=\sup\{t: E_Q(v^2_t)<\infty\} = \infty$. As a consequence, 
\begin{equation*}
E_Q(P(T)) < \infty \qquad \text{and} \qquad E_Q(P^n(T)) < \infty, 
\end{equation*}
and we have $\lim_{n\to\infty} P^n(T) = P(T)$.
\end{proposition}

\begin{proof}
We can find a Brownian motion $W^{\prime }$, independent of $W$ and $B$,
such that 
\begin{equation*}
Z_t = \rho W_t + \sqrt{1-\rho^2}W^{\prime }_t. 
\end{equation*}
Let $Q$ be the measure whose density process $Y_t=E_P\big(\frac{dQ}{dP}\mid%
\mathcal{F}_t\big)$ is given by $dY_t = -Y_t\gamma \sqrt{w_t}dW^{\prime }_t$%
, where $\gamma>0$ is a constant to be determined. To show that $Y$ is
indeed a martingale on $[0,T]$, it suffices to verify, as in Lemma~\ref%
{L:vmg}, that the auxiliary process 
\begin{equation*}
d\hat w_t=\big(\kappa(\theta-\hat w_t)-\gamma\rho\eta\hat w_t\big)dt+\eta%
\sqrt{\hat w_t}dZ_t, \qquad \hat w_0=w_0, 
\end{equation*}
is non-explosive. This can again be done using Feller's criteria. Next, it
follows from Girsanov's theorem that the dynamics of $w$ under $Q$ is given
by 
\begin{equation*}
dw_t = \tilde \kappa (\tilde \theta - w_t) dt + \eta\sqrt{w_t}d\tilde Z_t, 
\end{equation*}
where 
\begin{equation*}
\tilde \kappa = \kappa + \gamma\eta\sqrt{1-\rho^2}, \qquad \tilde\theta = 
\frac{\kappa\theta}{\kappa + \gamma\eta\sqrt{1-\rho^2}}, 
\end{equation*}
and $d\tilde Z_t = dZ_t + \gamma \sqrt{1-\rho^2} \sqrt{w_t} dt$ is Brownian
motion under $Q$. Hence, if we define 
\begin{equation*}
\tilde \chi = 2\rho\eta - \tilde\kappa \qquad \text{and} \qquad \tilde
\Delta = \tilde\chi^2 - 2\eta^2, 
\end{equation*}
we have that $\tilde T^*=\infty$ if $\tilde\chi\leq 0$ and $\tilde \Delta
\geq 0$. But 
\begin{equation*}
\tilde\chi = \chi - \gamma\eta\sqrt{1-\rho^2} 
\end{equation*}
and 
\begin{equation*}
\tilde\Delta = \Delta + \gamma\eta\sqrt{1-\rho^2}\big( \gamma\eta\sqrt{%
1-\rho^2} - 2\chi\big), 
\end{equation*}
so it suffices to choose $\gamma \geq \frac{2\chi}{\eta\sqrt{1-\rho^2}}$.

The verification of the last assertion is straightforward: $E_Q(P(T))<\infty$ is proved in the same way as under the measure $P$. To show that $E_Q(P^n(T))<\infty$ and $\lim_{n\to\infty} P^n(T) = P(T)$, note that $\int_0^T E_Q(v_t^2)dt <\infty$ due to the continuity and finiteness of $E_Q(v^2_t)$ on the compact interval $[0,T]$. An application of Theorem~\ref{T:c} concludes the proof.
\end{proof}




\subsection{The $3/2$-stochastic volatility model}

A model that has received considerable attention both in the theoretical and
empirical literature is the $3/2$-stochastic volatility process. See for
example \cite{Carr/Sun:2007} and the references therein. Let $B$ and $W$ be
two correlated Brownian motions. The model prescribes the following dynamics
for the stock price and its volatility 
\begin{align*}
dS_t &= S_t\sqrt{v_t}dB_t \\
dv_t &= v_t(p + q v_t)dt+\epsilon v_t^{\frac{3}{2}}dW_t
\end{align*}
where $p$, $q$ and $\epsilon$ are constants such that $q < \frac{\epsilon^2}{%
2}$ and $\epsilon>0$. The reason for the upper bound on $q$ is to avoid
explosion of $v$ in finite time. To see this, consider the process $R_t=%
\frac{1}{v_t}$, the reciprocal of $v$, which satisfies the SDE 
\begin{equation*}
dR_t=(\epsilon^2 - q - p R_t)dt-\epsilon \sqrt{R_t} dW_t. 
\end{equation*}
The process $R$ is a square-root process, and it is well-known that this
process avoids zero when $\epsilon^2-q>\frac{\epsilon^2}{q}$, which is
exactly the condition $q < \frac{\epsilon^2}{2}$. Let again $M$ be the
stochastic logarithm of $S$, i.e.~$M_t=\int_0^t\sqrt{v_s}dB_s$, so that $%
\langle M, M\rangle_t=\int_0^t v_sds$. Carr and Sun \cite{Carr/Sun:2007}
provide the Laplace transform of the integrated variance $\int_0^T v_sds$ in
closed form.

\begin{proposition}
\label{P:carr} In the $3/2$-model, the Laplace transform of the realized
variance $\int_0^T v_sds$ is given by 
\begin{equation*}
E(e^{-\lambda \int_0^Tv_sds})=\frac{\Gamma(\gamma-\alpha)}{\Gamma(\gamma)}%
\big(\frac{2}{\epsilon^2y_0}\big)^{\alpha}M(\alpha, \gamma, \frac{-2}{%
\epsilon^2y_0}) 
\end{equation*}
where $y_0=v_0\frac{e^{pT}-1}{p}$, $\alpha=-(\frac{1}{2}-\frac{q}{\epsilon^2}%
)+\sqrt{(\frac{1}{2}-\frac{q}{\epsilon^2})^2+2\frac{\lambda}{\epsilon^2}}$, $%
\gamma=2(\alpha+1-\frac{q}{\epsilon^2})$, $\Gamma$ is the Gamma function,
and $M$ is the confluent hypergeometric function 
\begin{equation*}
M(\alpha, \gamma, z)=\sum_{n=0}^{\infty}\frac{(\alpha)_n}{(\gamma)_n}\frac{%
z^n}{n!} 
\end{equation*}
with the notation $(x)_n=\prod_{i=0}^{n-1}(x+i)$.
\end{proposition}

\begin{proof}
We refer the reader to Carr and Sun \cite{Carr/Sun:2007}.
\end{proof}

Since the Laplace transform of the realized variance exists in a
neighborhood of zero, all moments of $\int_0^Tv_sds$ are finite. This
implies in particular that $E(\langle M, M\rangle_T) \in L^2$. From
Proposition \ref{P:fp}, both the true variance swap payoff and its
approximation have finite expectation.

Recall now the following result proved by Dufresne \cite{Dufresne:2001} on
the finiteness of moments of the square-root process.

\begin{proposition}
\label{P:duf} Let $\overline{v}=\frac{2(\epsilon^2-q)}{\epsilon^2}$. Then 
\begin{align*}
\forall \ p<\overline{v}, \qquad E(R_t^{-p}) < \infty \\
\forall \ p\geq\overline{v}, \qquad E(R_t^{-p}) = \infty
\end{align*}
and for all $p\geq -\overline{v}$, 
\begin{equation*}
E(R_t^p)=\mu_t^p e^{-\lambda_t}\frac{\Gamma(\overline{v}+p)}{\Gamma(%
\overline{v})}M(\overline{v}+p,p,\lambda_t) 
\end{equation*}
where $\mu_t=\frac{\epsilon^2}{2}\frac{1-e^{-pt}}{p}$, $\lambda_t=\frac{2pv_0%
}{\epsilon^2 (e^{-pt}-1)}$, $\Gamma$ is the Gamma function, and M is the
congruent hypergeometric function defined in Proposition~\ref{P:carr}.
\end{proposition}

If $q<0$, define $\kappa=-q>0$ and $\theta=\frac{p}{\kappa}$. Then the SDE
satisfied by $v$ can be re-written as 
\begin{equation*}
dv_t=\kappa v_t (\theta - v_t)dt + \epsilon v_t^{\frac{3}{2}}dW_t 
\end{equation*}
So under the condition $q<0$, the process $v$ is mean-reverting with a rate
of mean-reversion proportional to $v$. Also $\overline{v}>2$ when $q<0$, so
using Proposition~\ref{P:duf}, it can be seen that $E(v_t^2)=E(R_t^{-2})$ is
finite and integrable on $[0,T]$ as a continous function on this compact
time interval. Hence the condition of Theorem \ref{T:c} is satisfied and the
expectation of $P^n(T)$ converges to the expectation of $P(T)$ as $n
\rightarrow\infty$.

Under the condition $0\leq q< \frac{\epsilon^2}{2}$, it follows that $1<%
\overline{v}\leq 2$ and Proposition~\ref{P:duf} implies that $E(v_t^2)=\infty
$. By Fubini's theorem, $E(\int_0^Tv_s^2ds)=\infty$ so that the condition of
Theorem \ref{T:c} fails and the convergence of the $E(P^n(T))$ to $E(P(T))$
is not guaranteed anymore.

We now summarize the above findings.

\begin{example}
Suppose the stock price follows the $3/2$-stochastic volatility model. The
above discussion shows that

\begin{itemize}
\item[(i)] Both the true payoff $P^n(T)$ and the approximation $P(T)$ have
finite expectation.

\item[(ii)] If $q<0$, i.e.~when the squared volatility process is mean
reverting, $P^n(T)$ converges to $P(T)$ as $n\to\infty$.

\item[(iii)] If $q\geq 0$, our sufficient condition fails and we can no
longer guarantee that $P^n(T)$ converges to $P(T)$. It is an open problem to
establish whether or not this convergence actually takes place.
\end{itemize}
\end{example}


\end{document}